\documentclass[]{amsart}
\usepackage{bm}

\newcommand\defn[1]{{\bf #1}} 
\newcommand{\inv}{^{-1}}

\newcommand{\wZ}{\widetilde Z}
\newcommand{\wB}{\widetilde B}
\newcommand{\wC}{\widetilde C}

\newcommand{\bN}{\mathbb N}
\newcommand{\bZ}{\mathbb Z}
\newcommand{\bQ}{\mathbb Q}
\newcommand{\bR}{\mathbb R}
\newcommand{\bC}{\mathbb C}

\newcommand{\dra}{\dashrightarrow}

\newcommand{\dw}{{d}}
\DeclareMathOperator{\ndiff}{\operatorname{diff}}
\DeclareMathOperator{\Blowup}{\operatorname{Bl}}

\DeclareMathOperator{\codim}{\operatorname{codim}}
\DeclareMathOperator{\Exc}{\operatorname{Exc}}
\DeclareMathOperator{\Supp}{\operatorname{Supp}}
\DeclareMathOperator{\Sing}{\operatorname{Sing}}

\def\Div{\operatorname{Div}}

\newcommand\Q{{\mathbb{Q}}}
\newcommand\R{{\mathbb{R}}}

\theoremstyle{plain}
\newtheorem{theorem}{Theorem}[section]

\newtheorem{lemma}[theorem]{Lemma}
\newtheorem{corollary}[theorem]{Corollary}

\newtheorem{conditions}[theorem]{Conditions}
\theoremstyle{definition}

\newtheorem{definition}[theorem]{Definition}
\newtheorem{remark}[theorem]{Remark}

\newtheorem{example}[theorem]{Example}

\newtheorem*{acknowledgements}{Acknowledgments}

\author{Valery Alexeev}
\address{Department of Mathematics, University of Georgia, Athens GA
  30602, USA}

\author{Christopher Hacon}
\address{Department of Mathematics, University of Utah, Salt Lake
  City, UT 84112-0090, USA}

\author{Yujiro Kawamata}
\address{Department of Mathematical Sciences, University of Tokyo,
  Komaba, Meguro, Tokyo, 153-8914, Japan}  

\title{Termination of (many) 4-dimensional log flips}
\date{May 3, 2006}

\begin{document}
\begin{abstract}
  We prove that any sequence of 4-dimensional log flips that begins with
  a klt pair $(X,D)$ such that $-(K_X+D)$ is numerically equivalent to
  an effective divisor, terminates.  This implies termination of flips
  that begin with a log Fano pair and termination of flips in a
  relative birational setting.  We also prove termination of directed
  flips with big $K_X+D$.
  As a consequence, we prove existence of minimal models of
  4-dimensional dlt pairs of general type, existence of 5-dimensional
  log flips, and rationality of Kodaira energy in dimension 4.
\end{abstract}
\maketitle

Let $$(X,D) = (X^0,D^0)\dra (X^1,D^1)\dra (X^2,D^2)\dra \dots$$ be a
sequence of flips, relative over another variety $Z$, such that
\begin{enumerate}
\item $X^n$ are normal complex varieties of dimension 4, and $D^n=\sum
  d_jD_j^n$ are $\bR$-divisors,
\item $D^n$ are birational transforms of $D$ on $X^n$,
\item the pair $(X,D)$ is klt (Kawamata log terminal).
\end{enumerate}

The purpose of this note is to prove the termination of such sequences
in several situations. One of them is the case when $-(K_X+D)$ is
numerically equivalent over $Z$ to an effective $\bR$-divisor $E$,
i.e. $-(K_X+D)-E$ is an $\bR$-Cartier divisor which is zero on
curves contracted by $X\to Z$ (Theorem \ref{thm:negative-K}).  Note
that this effectivity condition holds in any of the following
situations:

\begin{enumerate}
\item $-(K_X+D)$ has nonnegative Kodaira dimension,
\item $X\to Z$ is generically finite. Indeed, in this case any
  divisor, including $-(K_X+D)$, is big over $Z$.
\end{enumerate}

Another situation is when $(X,D)$ is of general type over $Z$,
i.e. $K_X+D$ is big over $Z$. In this case we prove that a sequence of
\emph{directed} flips terminates (Theorem~\ref{thm:general-type}).
This allows us to prove that every 4-dimensional dlt pair of general type has a
minimal model (Corollary~\ref{cor:general-type}), and that
5-dimensional log flips exist (Theorem~\ref{5dim-flips}).
The paper is organized as follows:

In Section \ref{sec:preliminaries} we briefly review some
prerequisites.

Section~\ref{sec:Creative counting} is the heart of the paper. There,
we introduce a new \emph{difficulty} of a klt pair, improving on the
methods that were previously used by Shokurov, Kawamata, Matsuki,
Fujino and others \cite{Shokurov_Nonvanishing, KMM, Matsuki_Term4Flops,
  Kawamata_Term3logflips, Fujino_Termination, Kawamata_Term4Flips,
  Shokurov_Letters5}. Some key ideas for our definition come
from the unpublished manuscript \cite{Kawamata_Term4Flips} of 
Kawamata (using homology groups of the boundary divisors to store
information) and from the paper \cite{Shokurov_Letters5} of Shokurov
(using weights).  

A counting argument allows us to prove that for any fixed
$\alpha<1$ no divisors with discrepancy equal to $\alpha$ can appear
over flipping or flipped locus infinitely many times.
Then in Section~\ref{sec:Termination} we derive from this statement
the termination of flips in various useful situations. In the case
when $-(K_X+D)$ is effective, the termination is achieved in one
step. When $K_X+D$ is big, we first prove that sequences of
directed flips terminate. Then we establish the existence of a minimal
model. 

Section~\ref{sec:Applications} lists several consequences that follow
from our results. These include existence of 5-dimensional log flips
and the rationality of Kodaira energy for 4-dimensional klt pairs.

\begin{acknowledgements}
  A large part of this paper was written during the first author's
  2004 visit to RIMS, Kyoto University. He would like to express his
  gratitude for the hospitality and a very stimulating atmosphere. He
  was also partially supported by NSF grant 0401795.
  The second author was partially supported by NSF grant
  0456363 and by an AMS centennial scholarship.  
\end{acknowledgements}

\section{Preliminaries}
\label{sec:preliminaries}

We work over $\bC$ and use standard facts and notations of the Minimal
Model Program, as in \cite{KollarMori_Book}. Let us briefly recall
some of them.

Let $(X,D)$ be a pair consisting of a normal variety $X$ and an
$\bR$-Weil divisor $D=\sum d_jD_j$. It will be a standing assumption
that the divisor $K_X+D$ is $\bR$-Cartier. Then for every divisorial
valuation $v$ of the function field of $X$ one has the
\defn{discrepancy} (not the log discrepancy!)  $a(v, X,D)$ or simply
$a_v$ with respect to $K_X+D$.  The \defn{center of valuation}
$v$ on a variety $W$ will be denoted by $C(v,W)$. If $F=C(v,W)$ is a
divisor then we may write $a(F)$ for $a_v$.  For the components $D_j$
of $D$ one has $a(D_j, X,D) = -d_j$.

A pair $(X,D=\sum d_jD_j)$ is
\begin{enumerate}
\item \defn{terminal} if $d_j<1$ and for all $v\ne v_{D_j}$ one has
$a(v,X,D)>0$; we say that the
  pair is \defn{effective} if $d_j\ge 0$.
\item \defn{klt} (Kawamata log terminal) if $0<d_j<1$ and all $a(v,X,D)>-1$.
\end{enumerate}

A \defn{flip}, relative over $Z$, is a diagram
$
X^- \overset{\psi^-}{\longrightarrow} W
\overset{\psi^+}{\longleftarrow} X^+
$
of birational transformations of varieties over $Z$ such that
\begin{enumerate}
\item $\psi^-$ and $\psi^+$ are proper birational morphisms which are
  isomorphisms outside of closed subsets of codimension $\ge2$,
\item the relative Picard ranks of $X^-$ and $X^+$ over $W$ equal one,
\item $-(K_{X^-} + D)$ is $\psi^-$-ample and $K_{X^+} + D$ is $\psi^+$-ample.
\end{enumerate}
(Whenever two varieties are isomorphic in
codimension 1, we may denote birational transforms of divisors by the
same letters, for convenience.)
If $X^-$ is $\bQ$-factorial then so is $X^+$.  The exceptional set
$\Exc\psi^-$ is called the \defn{flipping locus} and $\Exc\psi^+$ the
\defn{flipped locus}.  

The most basic and well-known fact about flips is this:
\begin{lemma}\label{lem:flip-improves}
  Let $ X^- \overset{\psi^-}{\longrightarrow} W
  \overset{\psi^+}{\longleftarrow} X^+ $ be a flip. Then for any
  divisorial valuation $v$ one has $a(v,X^-,D^-) \le a(v,X^+,D^+)$ and
  the inequality is strict precisely for $v$ with $C(v,X^-) \subset
  \Exc\psi^-$, equivalently with $C(v,X^+)\subset \Exc\psi^+$.
\end{lemma}
In particular if $(X^-,D^-)$ is klt then so is $(X^+,D^+)$.

\begin{lemma}[\cite{KMM}, Lemma 5.1.17]
  \label{lem:types-of-flips}
  $\dim\Exc\psi^- + \dim\Exc\psi^+ \ge \dim X^- -1$. 
In particular, for $\dim X=4$ the possibilities for the pair 
$(\dim\Exc\psi^- , \dim\Exc\psi^+)$ are $(2,1)$, $(2,2)$ and $(1,2)$.
\end{lemma}

A \defn{terminalization} of a klt pair $(X,D)$ is a terminal pair
$(Y,B)$ together with a morphism $f:Y\to X$ such that $K_Y+B =
f^*(K_X+D)$. For any $v$ one has $a(v, X,D) = a(v, Y,B)$. All the
terminalizations used in this paper will be $\bQ$-factorial varieties.

\begin{lemma}\label{lem:terminalization}
  A $\bQ$-factorial terminalization exists for any klt pair $(X,D)$
  with $\dim X=4$.  
\end{lemma}
\begin{proof}
  For the proof, see \cite[Thm.5]{Kawamata_Term3logflips}. The proof
  is valid in any dimension provided the existence of flips (which for
  4-dimensional flips was proved by Shokurov and Hacon-McKernan
  \cite{Sho03,HaconMcKernan_Flips}) and the termination of terminal flips
  (Fujino \cite{Fujino_Termination}).
%
\end{proof}

A klt pair $(X,D)$ has only finitely many divisorial valuations with
$a_v\in (-1,0]$. They are:
\begin{enumerate}
\item the components $D_j$ of $D$, and 
\item the exceptional divisors of a terminalization $(Y,B)\to (X,D)$.  Let
$\bold{e(X,D)}$ denote the \defn{number of these exceptional
  divisors}.
\end{enumerate}
By Lemma~\ref{lem:flip-improves} in a sequence of flips the number
$e(X^n,D)$ 
decreases, and is eventually constant. 
There are more valuations in the next interval $(0,1)$:

\begin{example}\label{ex:echos}
  Let $(Y,B=\sum b_iB_i)$ be a terminal pair, and let $C\subset Y$ be a
  codimension-2 subvariety lying on a unique irreducible component
  $B_i$ of $B$ and not contained in $\Sing B_i$.
  Since $Y$ is terminal, it is nonsingular along the
  generic point of $C$.
  Let $E_1$ be the unique irreducible divisor of the blowup 
  $Y_1=\Blowup_{C}\overset{\pi}{\longrightarrow} Y$
  which dominates $C$, and let $C_1= E_1 \cap \pi_*\inv B_i$, the
  latter denoting the strict preimage of $B_i$. Let $E_2$ be the
  unique irreducible divisor of the blowup
  $Y_2=\Blowup_{C_1}Y_1$, etc. Then a simple computation (which is the
  same as in the case of a surface $Y$ and a point $C$) gives
  \begin{displaymath}
    a(E_1)= 1-b_i, \ a(E_2)= 2(1-b_i), \dotsc , \
    a(E_k)= k(1-b_i)
  \end{displaymath}
  We will call the divisorial valuation corresponding to $E_k$ the
  \defn{$\bold{k}$-th echo} of $Y$ along~$C$. 
\end{example}

\begin{lemma}\label{lem:echos}
  Let $(Y,B)$ be a terminal pair. Then the only divisorial valuations
  with $a(v, Y,B)\in (0,1)$ are:
  \begin{enumerate}
  \item The echos of $Y$ along irreducible subvarieties $C$ of
    codimension 2 that lie on a unique $B_j$, as above. For each such
    $C$, only finitely many $k$-th echos satisfying $k(1-b_j)<1$
    appear.
  \item Finitely many others.
  \end{enumerate}
\end{lemma}
\begin{proof}
  The proof is an easy application of the formula (see, e.g.,
  \cite{Kollar_Flops}) computing the discrepancies of a log smooth
  pair, applied to a log resolution of $(Y,B)$ on which the finitely
  many divisors $F_j$ with $a(Y,B,F_j)\le 0$ are disjoint.
\end{proof}


\begin{lemma}\label{lem:exc-divs-independent}
  Let $\pi:X\to Y$ be a birational morphism of projective varieties 
  with smooth $X$ and normal $Y$, and let
  $E_1,\dotsc,E_m$ be the irreducible components of the exceptional
  locus with $\codim_X E_i =1$. Then
  \begin{enumerate}
  \item $E_i$ are linearly independent in the Neron-Severi group $N^1(X)$
  \item The integer
    $$
    {\rho(Y)} = \dim N^1(X)\otimes\bR -
    \#(\text{of exceptional divisors } E_i)
    $$
    is well defined, i.e. does not depend on a particular resolution
    of $Y$, and is positive. 
  \end{enumerate}
\end{lemma}
\begin{proof}
  (1) Suppose that there is a nontrivial linear relation, and let $e$
  be the maximal dimension of $\pi(E_i)$ that appear. By cutting with
  $e$ general hyperplanes on $Y$ and $\dim X - e+2$ hyperplanes on $X$
  we can reduce to the case when $X$ and $Y$ are surfaces. By
  \cite{Mumford_SurfaceSingularities} the intersection form on the
  exceptional curves $E_i$ is negative definite, so no nontrivial
  linear relation between $E_i$ is possible.

  (2) Let $\phi:X'\to X$ be another resolution of singularities, with
  exceptional divisors $F_i$.  Then we have 
  homomorphisms $\phi_*$ and $\phi^*$ between the Neron-Severi groups.
  One has $\phi_*\phi^* D_1=D_1$ and $\phi^*\phi_*D_2 = D_2 + \sum n_j
  F_j$ for some $n_j\in \bZ$.  Therefore, $\dim N^1(X') -
  \dim N^1(X) = \#(F_j)$, and so $\rho(Y)$ is well defined. It is
  positive because $E_i$ are linearly independent and do not span
  $N^1(X,\bR)$. 
\end{proof}


\begin{remark}
  Our final remark is about dlt (divisorially log terminal) flips. A
  well-known statement called Special Termination (see e.g.
  \cite[2.1]{Fuj06}) reduces termination of dlt flips to the klt
  flips.  So all our termination results can be strengthened from klt
  to dlt.
\end{remark}


\section{Creative counting}
\label{sec:Creative counting}

Ever since Shokurov gave the original argument in the 3-dimensional
terminal case \cite{Shokurov_Nonvanishing}, the main idea for proving
termination of flips has remained the same. One introduces an appropriate
\emph{difficulty} $d(X,D)$ which counts valuations with certain
discrepancies, and then shows that $d(X^-,D^-) > d(X^+,D^+)$ for a
flip. Then one can conclude the termination if the difficulty is, say,
a nonnegative integer, or under somewhat weaker assumptions.

The main approach in this paper will be the same, with several twists.
As Example~\ref{ex:echos} and
Lemma~\ref{lem:echos} show, there are infinitely many valuations in $(0,1)$,
but except for the echos, there are only finitely many of them. The
new ideas are these: 
\begin{enumerate}
\item Count discrepancies in $(0,1)$ with a normalizing summand, so that the
  echos contribute zero to the sum, thus making it finite.
\item Count them with weights (cf. \cite{Shokurov_Letters5}).
\item Use the Picard groups of the boundary divisors $B_i$ on a
  terminalization $(Y,B)$ of $(X,D)$ appropriately to ``store'' some
  information (cf.  \cite{Kawamata_Term4Flips}).
\end{enumerate}

We first define the \defn{difficulty of a terminal pair}. That will be
a real number $\dw(X,D) $ which will depend on a 
\defn{weight function} $w:(0,+\infty)\to \bR_{\ge0}$ on the
set of discrepancies. The weight function is assumed to satisfy the following:
\begin{conditions}\label{conditions}
  \begin{enumerate}
  \item $w(a)\ge 0$ for all $a$, and $w(a)=0$ for $a\ge1$. 
  \item $w$ is decreasing (when we say \defn{increasing} or
    \defn{decreasing}, we always mean \defn{not necessarily strictly}).
  \item \label{item4b}
    Whenever $0< \sum m_i b_i < 1$ for some $m_i\in \bZ_{\ge0}$,
    one has
    $$w(1-\sum m_i b_i) \ge \sum m_i w(1-b_i).$$
  \item For $i=1,\dotsc, m$, let  $\{b_i^n,\ n\in\bN\}$ be $m$
    decreasing sequences of real numbers in $(0,1)$. Then there exists
    $n_0$ such that the following holds:

    Whenever $k\ge 2$, $m_i\ge0$ and $n\ge n_0$ are integers such that
    \begin{displaymath}
      1 - \sum m_i b_i^n \in (0,1), \quad\text{resp.}\quad
      k(1-b_1^n) - \sum_{i>1} m_i b_i^n  \in (0,1) 
    \end{displaymath}
    then one has
    \begin{displaymath}
      \ndiff_1(n) \ge \ndiff_1(n+1)  \quad\text{resp.}\quad
      \ndiff_2(n) \ge \ndiff_2(n+1),
    \end{displaymath}
    where
    \begin{align*}
&    \ndiff_1(n) = w(1 - \sum m_i b_i^n) - \sum m_i w(1-b_i^n) 
    \quad \text{and} \\
&    \ndiff_2(n) = 
    w\big( k(1-b_1^n) -  \sum_{i>1} m_i b_i^n \big) - w \big( k(1-b_1^n) \big)
    \end{align*}
  \end{enumerate}
\end{conditions}

For example, we could take for $w$ a piecewise linear function with
decreasing absolute values of slopes.
The main two functions we will use in this papers will be
$w_{\alpha}^-$ and $w_{\alpha}^+$, defined as follows:
\begin{definition}\label{def:functions}
  Let $\alpha\in (0,1)$. Then
  \begin{enumerate}
  \item $w_{\alpha}^-(x)=1-x$ for $x\le \alpha$ and
    $w_{\alpha}^-(x)=0$ for $x> \alpha$, and
  \item $w_{\alpha}^+(x)=1-x$ for $x<\alpha$ and $w_{\alpha}^+(x)=0$
    for $x\ge \alpha$.
  \end{enumerate}
  The reader is invited to do an elementary exercise now and check
  that both $w_{\alpha}^-$ and $w_{\alpha}^+$ satisfy the above
  conditions. The integer $n_0$ is chosen so that for $n\ge n_0$ for
  every $i$ the sequence $\{b_i^n\}$ settles into one of the open
  intervals between the points
  \begin{displaymath}
    0, \ 1-\alpha, \ 1-\frac{\alpha}{2},\ \dotsc\
    ,1-\frac{\alpha}{k},\ \dotsc
  \end{displaymath}
  or into one of these points.
\end{definition}

We will define another, \defn{summed weight} function $W:(-\infty,1)\to
\bR_{\ge0}$ by the formula
$$
W(b) = \sum_{k=1}^{\infty} w\big( k(1-b) \big). $$ It easily follows
that there are only finitely many nonzero terms in this sum, that
$W$ is a nonnegative and increasing function, and that $W(b)=0$ for
$b\le 0$.  The meaning of $w$ and
$W$ and the necessity of the above conditions will be clear from the
proofs of Lemmas \ref{lem:nonnegative}, \ref{lem:monotonicity1},
\ref{lem:monotonicity2}.

Let $\nu:\coprod \wB_i \to \cup B_i$ be the normalization of the
divisor $\Supp B$. For any irreducible subvariety $C\subset \Supp B$
the preimage $\wC = \nu\inv(C)$ splits into a union of irreducible
components $\wC_{ij}\subset\wB_i$.

\begin{definition}\label{defn:diff-terminal}
  Let $(Y,B)$ be a terminal pair, \emph{not necessarily effective,}
  i.e. possibly with some $b_i<0$. Then we define
  \begin{eqnarray*}
   \dw(Y,B) = 
   && \sum_{a(B_j)\le 0} W(b_i) \rho( \wB_i ) + 
   \sum_{v; \ \codim C(v,Y)\ne 2 } 
   w( a_v ) + \\
   && +\sum_{ {\rm irr.} C \subset Y; \codim C =2}
   \left[
     \sum_{v; \ C(v,Y)=C } 
     w( a_v ) - \sum_{\wC_{ij}} W(b_i) 
   \right]
 \end{eqnarray*}
\end{definition}

\begin{remark}
  For $w(x)=\max(1-x,0)$ this definition is a version of Shokurov's
  ``stringy'' invariant $\rho^2$ in \cite{Shokurov_Letters5} except
  that we take into account higher ``echos'' (which is necessary).
\end{remark}

\begin{lemma}\label{lem:well-defined}
  $\dw(Y,B)$ is well defined, i.e. only finitely many summands are
  nonzero.
\end{lemma}
\begin{proof}
  Indeed, if $C$ is a codimension-2 subvariety of $Y$ lying on a
  unique component $B_i$ and not in $\Sing B_i$, the contribution
  coming from the echo divisors of $C$ is zero. And by
  Lemma~\ref{lem:echos} there are only finitely many other nonzero terms.
\end{proof}

\begin{lemma}\label{lem:nonnegative}
  Assume that $(Y,B)$ is an effective terminal pair. Then 
  $\dw(Y,B) \ge 0$.
\end{lemma}
\begin{proof}
  Clearly, $\sum W(b_i) \rho( \wB_i )\ge 0$ and $w(a_v)\ge0$, so we
  only need to prove that the contributions from codimension-2
  subvarieties are nonnegative. This is a computation that can be done
  assuming that $B_i$'s are curves on a nonsingular surface $Y$
  passing through a nonsingular point $C$ with multiplicities $m_i$.
  Let $E_1$ be the exceptional divisor of the blowup of $Y$ at $C$. We
  have $a(E_1)= 1-\sum m_ib_i$. Let us order the divisors so that
  $b_1\ge b_2 \ge\dotsc $.

  \emph{Case 1: all $b_i<1/2$.} Then $W(b_i)=w(1-b_i)$. 
  Let $m_i'$ be the number of
  irreducible components $\wC_{ij}$ lying over $C$; then 
  $m_i' \le m_i$. Therefore, we obtain
  \begin{align*}
    \sum_{v;\ C(v,Y)=C } 
    w( a_v ) - \sum_{\wC_{ij}} W(b_i) \ge
    & \ w(1-\sum m_ib_i) - \sum m'_i w(1-b_i) \\
    \ge & \  w(1-\sum m_ib_i) - \sum m_i w(1-b_i) 
  \end{align*}
  which is nonnegative by Condition~\ref{conditions}(3).

  \emph{Case 2: $b_1\ge 1/2$.} Then, since $(Y,B)$ is terminal, we
  must have $m_1=1$ and $b_i<1/2$ for $i\ge2$. Then in the negative
  direction we have additional contributions coming from $W(b_1)$, and
  they are $w( k(1-b_1))$ for all $k\ge 2$.

  On the other hand, let $E_2$ be the exceptional divisor of the
  blowup at the intersection point of $E_1$ and the strict
  preimage of $B_1$; $E_3$ be the exceptional divisor of the blowup
  at the intersection point of $E_2$ and the strict preimage of $B_1$, etc. 
  Continuing this way by induction, for each $k\ge 2$ we obtain a
  divisor $E_k$ with
  $$
  a(E_k) =k(1-b_1) -\sum_{i>1} m_i^k b_i
  \quad \text{for some } m_i^k\ge 0.
  $$
  Since the function $w$ is decreasing by Condition
  \ref{conditions}(2), we have
  \begin{displaymath}
    w\big( a(E_k) \big) 
     \ge  w\big( k(1-b_1)\big).
  \end{displaymath}
  So the positive additional terms outweigh the negative ones.
\end{proof}

\begin{lemma}\label{lem:pullback}
  Let $f:Y'\to Y$ be a proper birational morphism and $(Y',B')$,
  $(Y,B)$ be two terminal pairs such that $K_{Y'}+B' = f^*(K_{Y}+B)$.
  Then $\dw(Y',B') = \dw(Y,B)$.
\end{lemma}
\begin{proof}
  Let $\wB'_i \to \wB_i$ be the induced birational morphisms between
  the normalizations of the divisors with $a(v)\le 0$. The only
  difference between $\dw(Y',B')$ and $\dw(Y,B)$ is in the
  codimension-2 subvarieties $\wC'_{ij}\subset \wB'_i$ whose image on
  $\wB_i$ has codimension $\ge3$. They appear in the first and in the
  last terms, and they cancel out by Lemma~\ref{lem:exc-divs-independent}. 
\end{proof}

\begin{definition}
  The \defn{difficulty of a klt pair} $(X,D)$ is defined to be 
  $\dw(X,D) = \dw(Y,B)$, where $f:(Y,B)\to (X,D)$ is any terminal pair
  such that $K_Y+B = f^*(K_X+D)$. 
\end{definition}

\begin{lemma}\label{lem:nonnegative2}
  If $(X,D)$ is a klt pair then $\dw(X,D)\ge 0$. 
\end{lemma}
\begin{proof}
  Indeed, we can use for $(Y,B)$ a terminalization of $(X,D)$, and
  then the statement follows by Lemma~\ref{lem:nonnegative}.
\end{proof}

\begin{lemma}[Monotonicity 1] \label{lem:monotonicity1}
  Let $(Y^1,B^1)$ and $(Y^2,B^2)$ be two birational terminal pairs
  that have the same ``$a(v)<0$'' parts, i.e. 
  for every component $B_i^1$ with $b_i^1>0$ there exists a
  birational to it component $B_i^2$ with $b_i^2=b_i^1$, and vice
  versa. Assume that one has
  \begin{displaymath}
    (\varphi^1)^*(K_{Y^1}+B^1) \ge (\varphi^2)^*(K_{Y^2}+B^2) 
  \end{displaymath}
  on a common log resolution $\varphi^n: U \to Y^n$. Then
  $\dw(Y^1,B^1) \ge \dw(Y^2,B^2)$.
\end{lemma}
\begin{proof}
  By Lemma \ref{lem:pullback} we can replace both pairs by their
  pullbacks on $U$. Then the only difference is between $w(a_v)$, and
  it is nonnegative because the weight function is decreasing. 
\end{proof}

\begin{lemma}[Monotonicity 2]\label{lem:monotonicity2}
  For $i=1,\dotsc, m$ let $\{b_i^n,\ n\in\bN\}$ be $m$ decreasing
  sequences of real numbers in $(0,1)$, and let $n_0$ be chosen as in
  Condition~\ref{conditions}(4).  Let $(Y,B^n=\sum b_i^n B_i)$ be the
  effective terminal pairs with the same variety $Y$.  Then for $n\ge
  n_0$ one has $\dw(Y,B^n)\ge \dw(Y,B^{n+1})$.
\end{lemma}
\begin{proof}
  We clearly have the required inequality for the first term $\sum
  W(b_i)\rho(\wB_i)$. Also, for all valuations with $\codim_Y
  C(v,Y)>2$ one has $a(v,Y,B^n) \le a(v,Y, B^{n+1})$ and so for the
  corresponding terms in Definition \ref{defn:diff-terminal} one gets
  $w(a(v,Y,B^n)) \ge w(a(v,Y, B^{n+1}))$.

  It remains to consider valuations with $\codim_Y C(v,Y)=2$. As in
  the proof of Lemma~\ref{lem:nonnegative}, we have
  \begin{enumerate}
  \item the terms 
    \begin{align*}
      &    w(1 - \sum m_i b_i^n) - \sum m_i w(1-b_i^n) 
      \quad \text{and} \\
      &    
      w\big( k(1-b_1^n) -  \sum_{i>1} m_i^k b_i^n \big) 
      - w\big( k (1-b_1^n) \big),
    \end{align*}
  \item   and then some other positive terms of the form 
    $w(c - \sum m''_i b_i^n)$ for some $m''_i\ge 0$. (We used the fact
    that $m'_i \le m_i$.)
  \end{enumerate}
  For the terms in (1) we have the required inequality by
  Condition~\ref{conditions}(4). For the terms in (2) the inequality
  follows from the fact that $w$ is a decreasing function.
\end{proof}


\begin{theorem}\label{thm:klt-flips-decrease}
  Let $(X,D) = (X^0,D)\dra (X^1,D)\dra (X^2,D)\dra \dots$ be a
  sequence of flips with $(X,D)$ klt.  Then there exists $n_0$ such
  that for $n\ge n_0$ the sequence $\dw(X^n,D)$ is decreasing.
\end{theorem}
\begin{proof}
  Let $ X^n \overset{\psi^-}{\longrightarrow} W
  \overset{\psi^+}{\longleftarrow} X^{n+1} $ be the $n$-th flip.  For
  $n$ large enough, we have $e(X^n,D) = e(X^{n+1},D)$, i.e. on the
  terminalizations $(Y^n,B^n)$, $(Y^{n+1},B^{n+1})$ the irreducible
  components of $B^{n+1}$ are birational transforms of the irreducible
  components of $B^n$. Then we have sequences $\{b_i^n\}$ of
  coefficients of $B^n$, and we choose $n_0$ large enough, as in
  Lemma~\ref{lem:monotonicity2}. 

  Let $\varphi^-:U_n\to Y^n$, $\varphi^+:U_n\to Y^{n+1}$ be a variety
  dominating both terminalizations. Denoting the birational transform
  $(\psi^+_*)\inv \psi^-_* B^n$ on $Y^{n+1}$ again by $B^n$, we have
  $B^n \ge B^{n+1}$ on $Y^{n+1}$ by Lemma~\ref{lem:flip-improves}.
  Further, we claim that
  \begin{displaymath}
    (\varphi^-)^* (K_{Y^n} + B^n) \ge 
    (\varphi^+)^* (K_{Y^{n+1}} + B^n) 
  \end{displaymath}
  Indeed, the first divisor is numerically negative over $W$, and
  therefore also over $Y^{n+1}$, and the second divisor is zero over
  $Y^{n+1}$.  Hence, the difference is negative over $Y^{n+1}$, and it
  consists of exceptional divisors of $\varphi^+$. This implies that
  the difference is effective.

  Now, we have $\dw(Y^n, B^n) \ge \dw(Y^{n+1}, B^n)$ by
  Lemma~\ref{lem:monotonicity1} and $\dw(Y^{n+1}, B^n) \ge \dw(Y^{n+1},
  B^{n+1})$ by Lemma~\ref{lem:monotonicity2}.
\end{proof}

To summarize what we have done this far: For every weight function $w$
satisfying Conditions~\ref{conditions} and for any klt pair $(X,D)$ we
have defined difficulty $\dw(X,D)\ge 0$ and proved that in a sequence
of flips it decreases at each step for $n\ge n_0$. Now, if we could
only establish that it decreases by at least a fixed $\epsilon>0$, that would
certainly imply termination. 
So far we have not used the fact that for some valuations $v$, the
weight $w(a_v)$ may drop abruptly. We only used the fact that after a flip one
has $a_v(X^-,D^-) \le a_v(X^+,D^+)$ and so $w(a_v(X^-,D^-)) \ge
w(a_v(X^+,D^+))$.

\begin{theorem}\label{thm:fixed-A}
  Fix a number $\alpha\in (-1,1)$. Then in a sequence of flips 
  $$(X,D) = (X^0,D^0)\dra (X^1,D^1)\dra (X^2,D^2)\dra \dots$$ 
  there cannot be infinitely many $k$ for which there exist valuations
  $v$ with $a_v=\alpha$ whose center is in the flipping or the flipped
  locus.
\end{theorem}
\begin{proof}
  If $\alpha\in (-1,0]$, the situation is easy, since there are only
  finitely many discrepancies $a_v\in(-1,0]$ and after a flip $a_v$
  strictly increases. 
  
  So assume $\alpha\in (0,1)$.  For the flipping locus we use the
  function $w_{\alpha}^-$ of \ref{def:functions}. Then after a flip
  the corresponding discrepancy changes from $a_v(X^n,D^n)=\alpha$ to
  $a_v(X^{n+1},D^{n+1})>\alpha$.  Hence, the weight
  $w_{\alpha}^-(a_v)$ and the difficulty drop down by $(1-\alpha)$.
  This can happen only finitely many times. For the flipped locus we use
  the function $w_{\alpha}^+$ instead.
\end{proof}

\begin{lemma}\label{lem:no-sing-codim2}
  In a sequence of klt flips it cannot happen infinitely many times
  that the flipping or flipped locus has a component which has
  codimension $2$ in $X^n$ and is contained in the singular
  locus of $X^n$.
\end{lemma}
\begin{proof}
  Let $E_1, \dots, E_m$, $m=e(X^n,D)$, be the finitely many valuations
  of $X^n$ with $a(v,X^n,D)\le 0$ for $n\ge n_0$.

  Assume that $X^n$ is generically singular along a codimension-2
  component. By cutting with two generic hyperplanes we obtain a log
  terminal singular surface pair $(Z,D|_Z=\sum d_j D_j|_Z)$, where
  $D_j|_Z$ need no longer be irreducible. Let $\wZ$ be the minimal
  resolution of singularities, with exceptional curves $F_k$. Then for
  each $F_i$, $a(F_i,Z,D|_Z) \le 0$, and so $a(F_i,Z,D|_Z) = a(E_i,
  X^n,D)$ for some $E_i$. Since there are only finitely many $E_i$'s,
  there will be a subsequence such that for one of the divisors, say
  $E_1$, the discrepancy $a(E_1,X,D)$ appears as the minimal
  discrepancy of $(Z,D|_Z)$, and increases on every step. But the set
  of minimal discrepancies of log terminal pairs satisfies the
  ascending chain condition by \cite[Thm.3.2,3.8]{Alexeev_22Terms}.
  Contradiction.

  \emph{Second proof}. The positive numbers $-F_i^2$ are bounded from
  above by $2/(1+ \min(a_v(X,D)))$. If we could bound the number of
  $F_i$'s, there would be only finitely many possibilities for the
  weighted graphs of minimal resolutions of $Z$'s, and so finitely
  many possible indexes, and finitely many possible discrepancies
  $a(F_i)$. The statement then would easily follow.  

  Each divisor $E_i$ can lead to many $F_i$'s, so we cannot limit the
  number of curves on the minimal resolution of $Z$. However, it is
  still true that the number of \emph{distinct} discrepancies
  $a(F_i,Z,D|_Z)$ is $\le m$. Now, the dual graph of the minimal
  resolution of a singularity on $Z$ is a tree with at most one fork
  and three legs
  (by the classification of log terminal surface singularities, see
  e.g. \cite{Kawamata_Crepant,Alexeev_LogCanSings}).
  A basic computation shows that on the interior of each leg
  \begin{enumerate}
  \item the function $a(F_i)$ is concave up, and
  \item if $a(F_{i-1})=a(F_i)=a(F_{i+1})$ then $F_i^2=-2$ and $F_i\cdot D|_Z=0$.
  \end{enumerate}
  Hence, the only way one can have repeating discrepancies in the
  interior of a leg is when you have chains of $(-2)$-curves. These
  chains can be shortened without changing the set of
  discrepancies. Hence, we can assume that in the interior of each leg
  you have no more than 2 repeating $a(F_i)$'s. We are reduced to the
  finitely many graphs with at most $3m+7$ vertices, and so we are
  done by the above argument.
\end{proof}

\begin{theorem}\label{thm:no-codim2}
  In a sequence of klt flips it cannot happen infinitely many times
  that the flipping or flipped locus has a component which has
  codimension $2$ in $X^n$ and which is contained in $D^n$.
\end{theorem}
\begin{proof}
  Let $C$ be such a codimension-2 component.  By the previous lemma we
  can assume that $Y$ is nonsingular generically along $C$. The blowup
  along $C$ produces a divisor with $a_v = 1-\sum m_ib_i <1$. We have
  a finite set of such discrepancies in $(-1,1)$, and we apply
  Theorem~\ref{thm:fixed-A} to each of them to complete the proof.
\end{proof}

\section{Termination and minimal models}
\label{sec:Termination}

\begin{lemma}\label{lem:positive-bdry}
  Let $$(X,D) = (X^0,D)\dra (X^1,D)\dra (X^2,D)\dra \dots$$ be a
  sequence of 4-dimensional klt flips such that $D=E + D'$ with
  effective $E$ and $D'$, such that all flips are $E$-positive (i.e.
  $E$ is positive on a curve in the flipping locus).  Then the
  sequence of flips terminates.
\end{lemma}
\begin{proof}
  Write $E=\sum e_jD_j$ as the sum of irreducible components, with
  $e_j>0$. Then for every flip there exists a component $D_j$ such
  that $D_j$ is positive on the flipping locus, and so is negative on
  the flipped locus. So, $D_j$ must contain the flipped locus. 

  By Theorem~\ref{thm:no-codim2} after some $n\ge n_0$ the flipped
  locus does not contain any codimension-2 components. By
  Lemma~\ref{lem:types-of-flips} this implies that all the flips for
  $n\ge n_0$ are the $(2,1)$-flips. But for each of these flips the
  dimension of the homology group $H_4(X^n,\bR)$ goes down, and this
  cannot be repeated infinitely many times.
\end{proof}

\begin{theorem}\label{thm:negative-K}
  Let $$(X,D) = (X^0,D)\dra (X^1,D)\dra (X^2,D)\dra \dots$$ be a
  sequence of klt flips, relative over another variety $Z$.
  Assume that $-(K_X+D)$ is numerically equivalent over $Z$ to an effective
  $\bR$-divisor $E$ 
  (i.e. the difference is an $\bR$-Cartier divisor which is
  zero on curves contracted by $X\to Z$).
   Then the sequence of flips is finite.
\end{theorem}
\begin{proof}
  Note that $E\ne0$. 
  Replace $D$ by $D+ \epsilon E$. Then for $0<\epsilon \ll1$ the pair
  $(X,D+\epsilon E)$ is still klt and since $K+D+\epsilon E \equiv
  (1-\epsilon)(K+D)$, the sequence of flips for $K+D$ is a sequence of
  flips for $K+D+\epsilon E$. Since $E$ is positive on each flipping
  locus, we are done by the previous Lemma~\ref{lem:positive-bdry}.
\end{proof}

Our next application is to the case when $K_X+D$ or $D$ are big.  Let
us begin by recalling the following consequence of the Cone Theorem
(which follows immediately from \cite{KMM}):

\begin{lemma}\label{L1} 
  Let $(X,D)$ be a klt pair and $f:X\to Z$ a projective morphism to a
  normal variety.  Let $A$ be an effective $\R$-Cartier divisor such that
  $K_X+D$ is not relatively nef, but $K_X+D+A$ is relatively nef. Then
  there is a $K_X+D$ extremal ray $R$ over $Z$ and a real number
  $0<\lambda \leq 1$ such that $K_X+D+\lambda A$ is relatively nef,
  but trivial on $R$.

  If $D$ and $A$ are $\bQ$-Cartier divisors then $\lambda\in \bQ$.
\end{lemma}

Given $(X,D)$ as above and an ample line bundle $A$ on $X$, one may
run the \defn{MMP for $\bold{(X,D)}$ directed by $\bold{A}$.} 
Recall here that 4-dimensional log flips exist by Shokurov and 
Hacon-McKernan \cite{Sho03,HaconMcKernan_Flips}.
We proceed as follows: Since
$A$ is ample, then $K_X+D+tA$ is nef for some $t\gg 0$. Let
\begin{displaymath}
  t_0=\inf \{ t\geq 0 | K_X+D+tA\ {\rm is\  relatively\ nef }\}.
\end{displaymath}
If $t_0=0$, or $R$ induces a Mori-Fano fibration, we stop,
otherwise we replace $(X,D)$ by the pair $(X^1, D)$ given by flipping
or contracting $R$.  Now $K_{X^1}+D+t_0A$ is relatively nef and by
\eqref{L1} we can repeat the above trick.  At each step, we obtain a
pair $(X^n, D)$ and a real number $t_n$ such that $K_{X^n}+D+t_{n}A$
is relatively nef and $t_n$ is the infimum of all $t\geq 0$ such that
$K_{X^n}+D+tA$ is relatively nef.  Note that $X^{n+1}$ is obtained
from $X^n$ by performing a $(K_{X^n}+D)$-flip or divisorial contraction
with extremal ray $R_n$.

It is clear that each flip or divisorial contraction is $(K+D)$-negative
and $A$-positive and we have the following possibilities:
\begin{enumerate}
\item $K_{X^n}+D$ is relatively nef;
\item $R^n$ induces a relative Mori-Fano fibration for $K_{X^n}+D+t_nA$ and
  hence for $K_{X^n}+D$;
\item For all $n\geq n_0$, $R^n$ induces a $(K_{X^n}+D)$-flip.
\end{enumerate}
  
\begin{theorem}\label{thm:general-type}
  Let $(X,D)$ be a 4-dimensional klt pair over $Z$.
  Assume that for some $c_j\in \bR$ the $\bR$-divisor $c_0 K_X + \sum
  c_jD_j$ is numerically equivalent to a big divisor.
  Then there exists a finite sequence of flips and divisorial
  contractions $X^n\dasharrow X^{n+1}$, $0\leq k\leq n$
  (here $X^0=X$) such that $K_{X^n}+D$ is nef 
  or there exists a $K_{X^n}+D$ Mori-Fano fibration.
\end{theorem}
\begin{proof}
  Denote $D'=\sum c_jD_j$. Write the big divisor as $A+E$ with an
  effective divisor $E$ and an ample divisor $A$ whose support does
  not have components in common with $D+E$. Then
    \begin{eqnarray*}
&    K+D \equiv 
    (1-\epsilon c_0)K + D-\epsilon D' + \epsilon A + \epsilon E
    =
    (1-\epsilon c_0)(K + D''), \\
&    \text{where}\quad
    D'' = (1-\epsilon c_0)\inv (D-\epsilon D' + \epsilon A + \epsilon E)
    \end{eqnarray*}

For $0<\epsilon\ll 1$ the pair $(X,D'')$ is klt, the divisor $D''$ is
effective and contains a positive multiple of $A$. If we run the MMP
directed by $A$ then each extremal ray will be $A$-positive. So by
Lemma~\ref{lem:positive-bdry} the flips terminate. So we either arrive
to a minimal model of $K+D''$, and hence of $K+D$ which is
proportional to it, or we end up with a Mori-Fano fibration.
\end{proof}

\begin{corollary}\label{cor:dltgt}
Let $(X,D)$  be a $4$-dimensional $\Q$-factorial dlt pair, such that
$D\in \Div _\Q(X)$ is big. Let $A$ be any ample divisor on $X$.
Then the MMP for $(X,D)$ directed by $A$ terminates.
\end{corollary}
\begin{proof} It suffices to show that there is no infinite sequence of 
$K_X+D$ flips (directed by $A$) $(X^n,D)\dasharrow (X^{n+1},D)$.
By \cite[2.1]{Fuj06}, we know that after finitely many flips, the flipping 
and flipped loci are disjoint from the locus of log canonical singularities 
$\lfloor D\rfloor $. We now conclude as in the proof of Theorem
\ref{thm:general-type}. 
\end{proof}

\begin{corollary}\label{cor:general-type}
  Let $(X,D)$ be a dlt pair of general type over $Z$. Then there exists
  a minimal model of $(X,D)$ over $Z$.
\end{corollary}
\begin{proof}
  Indeed, the Mori-Fano fibration in this case is not possible. 
\end{proof}

\section{Further applications}
\label{sec:Applications}

We begin by recalling a result concerning the big cone of a
klt pair.
Let $(X,D)$ be a projective $\Q$-factorial klt pair, $D\in \Div_\Q (X)$ and $A$ an ample divisor on $X$.
Then one has:
\begin{definition} If $K_X+D$ is not pseudo-effective, then 
the \defn{effective log-threshold} of $A$ is defined by
$$\sigma (X,D,A)={\operatorname {sup}}\{t\in \Q \ |\ A+t(K_X+D)\ {\rm is\ effective}\}.$$ The \defn{Kodaira Energy} is defined by $\kappa\epsilon (X,D,A)=-1/\sigma (X,D,A)$.\end{definition}
Following the approach of Batyrev \cite{Batyrev_EffDiv}, we deduce the
following result concerning the rationality of the Kodaira Energy (see
\cite{Ara05} for more details).
\begin{theorem}\label{thm:KE} Let $(X,D)$ be a projective $4$-dimensional
$\Q$-factorial klt pair and $A$ an ample divisor on $X$.
If $K_X+D$ is not pseudo-effective, then $\sigma (X, D , A)$ is rational.
\end{theorem}
\begin{proof} 
  We use the MMP for $(X,D)$ directed by $A$.  Since $K_X+D$ is not
  pseudo-effective, it ends with a Mori fiber space.  Let $R^n$ be the
  final extremal ray on the model $(X^n,D)$.  Let $\lambda$ be the
  nef threshold which is known to be rational.  Then
  $K_{X^n}+D+\lambda A^n$ is nef, hence pseudo-effective, and
  $(K_{X^n}+D+\lambda A^n) \cdot R^n=0$.  Since $-(K_{X^n}+D)$ is
  relatively ample for the Mori fiber space, we conclude that 
  $\sigma(X,D,A)=-1/\lambda\in \bQ$ as required.
%
%

\end{proof}
We now turn our attention to the problem of the existence of flips in
dimension~$5$. We have the following:

\begin{theorem}\label{5dim-flips}
  Let $(X,D)$ be a $5$-dimensional klt $\Q$-factorial pair, $D\in \Div
  _\Q (X)$ (or $\Div_\R(X)$), $f:(X,D)\to Z$ a flipping contraction.
  Then the flip $f^+:X^+\to Z$ of $f$ exists.
\end{theorem}
In \cite{HaconMcKernan_Flips}, the existence of flips in dimension $n$
is claimed assuming the termination of flips in dimension $n-1$.
Since we have not shown the termination of flips for $4$-folds in full
generality, we must analyze the given argument in more detail. 
\begin{proof} By \cite{HaconMcKernan_Flips}, PL-flips exist in dimension $5$,
provided that flips (relative over a variety $Z$) 
terminate for any klt pair $(T, B)$ where $T\to Z$ is birational and
$\dim T=4$. This is the case by Theorem \ref{thm:negative-K}.
We now follow Shokurov's Reduction Theorem as outlined in \S 3 of 
\cite{Fuj06}.
We do not reproduce the entire argument, but simply indicate the changes 
required to 3.7 of \cite{Fuj06}.

In Step 2 of 3.7, we have
$\nu:Y\to X$ an appropriate log resolution of $X$ and 
a morphism $h:Y\to X\to Z$. 
$H'$ is an appropriately chosen divisor on $Z$ and $H:=f^*H'$.
We may assume that $D$ contains a general $f$-ample $\Q$-divisor $D'$. We let
$A\sim _{\Q} e_0\nu ^*D'+\sum e_iE_i$ be a general ample divisor over $Z$ 
where $0<e_0\ll 1$, $0<-e_i\ll 1$ and $E_i$ are
$\nu$-exceptional divisors. 
Note that $K_Y+(D+H)_Y+\sum e_iE_i\sim_\Q K_Y +(D^0+H)_Y+A$ where $D^0=D-e_0D'\geq 0$ and 
$(D+H)_Y$ denotes the birational transform of $D+H$, 
that is the strict transform of $D+H$
plus the sum of all the exceptional divisors taken with multiplicity $1$.
Notice that we may replace $D$ by $D-e_0D'+\nu _*A$ so that
$K_Y +(D^0+H)_Y+A-\nu ^*(K_X+D)$ is given by the strict transform of $H'$ plus an exceptional divisor.

In Step 3 we run the MMP with respect to $K_Y+(D^0+H)_Y+A$ over $Z$ directed by $A$. All flips are PL-flips and their
termination follows by the usual special termination argument (as in \S 2 of \cite{Fuj06})
noting that to run this argument, it suffices to establish termination of directed flips for dlt pairs with big boundary divisor
in dimension $\leq 4$. This in turn follows from 
Corollary \ref{cor:dltgt}.
We obtain a $\Q$-factorial dlt pair $\bar{h}:(\bar{Y},(D^0+H+A)_{\bar {Y}})\to Z$ such that $K_{\bar{Y}}+(D^0+H+A)_{\bar{Y}}$ is $\bar{h}$-nef.

In Step 5, we subtract $\hat {H}$ (where $\hat{H}$ denotes
the strict transform of $H$ on $\bar{Y}$).
The main point is that for this step to work, we need the special termination
of dlt
flips directed by $\hat {H}$. By the arguments of Special Termination (\S 2 of
\cite{Fuj06}), we need the termination of dlt
flips directed by the restriction $\hat {H}$ for log pairs of dimension
at most $4$ whose boundary is a big $\Q$-divisor. 
This follows as in Corollary \ref{cor:dltgt}.
Therefore we obtain $\tilde{h}:(\tilde{Y},D_{\tilde {Y}}^0+\tilde{A})\to Z$ such that $\tilde {Y}$ is
$\Q$-factorial, $(\tilde{Y},D_{\tilde{Y}}^0+\tilde{A})$ is dlt and $K_{\tilde{Y}}
+D_{\tilde{Y}}^0+\tilde{A}$ is $\tilde{h}$-nef.
By the negativity lemma, one can check that $\tilde{h}$ is small
so that $(\tilde{Y},D_{\tilde{Y}}^0+\tilde{A})$ is in fact klt.
\footnote{To see this, let $W$ be a common log resolution of $\tilde{Y}$ and $X$ with morphisms $\phi :W\to \tilde{Y}$ and $\psi:W\to X$. Then 
$-B:=\phi ^*(K_{\tilde{Y}}+D^0_{\tilde {Y}}+\tilde{A})-\psi ^*(K_X+D)$ is exceptional and $\psi$-nef, therefore $B$ is effective and $\phi ^*(K_{\tilde{Y}}+D^0_{\tilde {Y}}+\tilde{A})+B=\psi ^*(K_X+D)$. If $E\subset \tilde{Y}$ is a divisor which is $Z$-exceptional, then this implies that $a(E,X,D)\leq -1$ which is absurd. It also follows that the canonical models over $Z$ for $(X,D)$ and $(\tilde{Y},D^0_{\tilde {Y}}+\tilde{A})$ coincide.}
One sees that the canonical model of 
$(\tilde{Y},D_{\tilde{Y}}^0+\tilde{A})$ over $Z$ coincides with the
canonical model of $(X,D)$ over $Z$ and this is 
the required flip.
\end{proof}

\bibliographystyle{amsalpha}

\def\cprime{$'$}
\providecommand{\bysame}{\leavevmode\hbox to3em{\hrulefill}\thinspace}
\providecommand{\MR}{\relax\ifhmode\unskip\space\fi MR }
\providecommand{\MRhref}[2]{%
  \href{http://www.ams.org/mathscinet-getitem?mr=#1}{#2}
}
\providecommand{\href}[2]{#2}

\end{document}